\documentclass[12pt]{amsart}

\usepackage{amssymb}
\usepackage{amsmath}
\usepackage{amsthm}
\usepackage{a4wide}

\newtheorem{theorem}{Theorem}

\newtheorem{proposition}[theorem]{Proposition}

\newtheorem{conjecture}[theorem]{Conjecture}

\newtheorem{question}[theorem]{Question}

\def\leq{\leqslant}
\def\geq{\geqslant}
\def\al{\alpha}
\def\be{\beta}

\def\la{\lambda}
\def\C{\mathbb C}

\def\N{\mathbb N}

\def\Q{\mathbb Q}
\def\Z{\mathbb Z}

\def\Om{\Omega}

\title[On the equation $f(g(x))=f(x)h^m(x)$]{On the equation $f(g(x)) =f(x)h^m(x)$ for composite polynomials}

\author{Himadri Ganguli}
\address{Department of Mathematics\\Simon Fraser University\\8888 University Drive\\Burnaby, British Columbia V5A 1S6\\Canada}
\email{hganguli@sfu.ca}

\author{ Jonas Jankauskas}
\address{Department of Mathematics and Informatics\\ Vilnius University\\
Naugarduko 24, Vilnius LT-03225\\ Lithuania}
\email{jonas.jankauskas@gmail.com}
\thanks{A visit of the second author at IRMACS Center, Simon Fraser University was funded by Lithuanian Research Council (Student research support project).}

\subjclass[2000]{11B83, 11C08, 11D57, 11N32, 11R09, 12D05, 12E10}

\keywords{Chebyshev polynomial, composite polynomials, Pell equation, multiplicative dependence}
\begin{document}

\maketitle

\begin{abstract} In this paper we solve the equation $f(g(x))=f(x)h^m(x)$ where $f(x)$, $g(x)$ and $h(x)$ are unknown polynomials with coefficients in an arbitrary field $K$, $f(x)$ is non-constant and separable, $\deg g \geq 2$, the polynomial $g(x)$ has non-zero derivative $g'(x) \ne 0$ in $K[x]$  and the integer $m \geq 2$ is not divisible by the characteristic of the field $K$. We prove that this equation has no solutions if $\deg f \geq 3$. If $\deg f = 2$, we prove that $m = 2$ and give all solutions explicitly in terms of Chebyshev polynomials. The diophantine applications for such polynomials $f(x)$, $g(x)$, $h(x)$ with coefficients in $\Q$ or $\Z$ are considered in the context of the conjecture of Cassaign et. al on the values of Louiville's $\lambda$ function at points $f(r)$, $r \in \Q$.
\end{abstract}

\section{Introduction}\label{intr}

The problem investigated in the present paper is motivated by the following question:
\begin{question}\label{klaus}
Do there exist integer polynomials $f(x)$, $g(x)$ and $h(x)$ of degrees $\deg f \geq 3$, $\deg g \geq 2$, $f(x)$ separable (and possibly irreducible in $\Z[x]$), such that $f(g(x))=f(x)h^2(x)$?
\end{question}
\noindent This question has been posed in connection with a recent work of Borwein, Choi and Ganguli \cite{lui} on the sign changes of the \emph{Liouville's lambda function} $\la(f(n))$ for the values of integer quadratic polynomials $f(x) \in \Z[x]$ at integer points $n \in \Z$. Recall that for $n \in \Z$, the lambda function $\la(n)$ is defined by $ \la(n) = (-1)^{\Om(n)}$, where $\Om(n)$ is the total number of prime factors of $n$, counted with multiplicity. Alternatively, $\la(n)$ is the completely multiplicative function defined by $\lambda (p)=-1$ for each prime $p$ dividing $n$. Chowla \cite{chow} conjectured that \[\sum_{n \leq x}\la(f(n)) = o(x)\] for any integer polynomial $f(x)$ which is not of the form $f(x) = b g(x)^2$, where $b \in \Z$ and $g(x) \in \Z[x]$. For $f(x) = x$, Chowla's conjecture is equivalent to the prime number theorem and has been proven for linear polynomials $f(x)$, but is open for polynomials of higher degrees. Even the much weaker conjecture of Cassaigne et al. \cite{cass} which states
\begin{conjecture}
If $f(x) \in \Z [x]$ and is not of the form of $bg^2(x)$ for some $g(x)\in \Z[x]$, then $\lambda (f(n))$ changes sign infinitely often.  \end{conjecture}\noindent has not been proved unconditionally for the polynomials of degree $\deg f \geq 2$.

In the paper \cite{lui} it has been proved that the sequence $\la(f(n))$ cannot be eventually constant for quadratic integer polynomials $f(x)=ax^2+bx+c$, provided that at least one sign change occurs for $n >(|b|+(|D|+1)/2)/2a$, where $D$ is the discriminant of $f(x)$. The proof is based on the solutions of Pell-type equations. In practice, using this conditional result, one can prove the Cassaigne's conjecture for any particular integer quadratic $f(x)$, for instance, $f(x)=3x^2+2x+1$. In contrast, the only examples of degree $\deg f \geq 3$ for which the conjecture has been proven in \cite{cass}  are $f(x) =\prod_{j=1}^{k}(ax+b_j)$, where $a, b_k \in \N$, $b_k$ are all distinct, $b_1 \equiv \dots \equiv b_k \pmod{a}$. No similar examples of irreducible integer polynomials of degree $d \geq 3$ are known . It appears that the problem of finding an irreducible example of degree $d=3$ is interesting and probably difficult. 

We now explain how the composition identity in Question \ref{klaus} could be of use to prove that $\la(f(n))$ or $\la(f(-n))$ is not eventually constant for cubic polynomials $f(x)$. Assume that the leading coefficient of $g(x)$ is positive. Since $\deg g \geq 2$, there exists a positive integer $n_0$ such that $g(n)>n$ for integers $n > n_0$. Suppose that there exist two integers $k_0, l_0 > n_0$ such that $\la(f(k_0))=-\la(f(l_0))$. Then $\la(f(k_j))$ and $\la(f(l_j))$ also differ in sign for infinite sequences  of integers $k_j$ and $l_j$, defined by $k_{j+1} = g(k_j)$ and $l_{j+1}=g(l_{j})$, $j \geq 0$, since $\la(f(g(n)))=\la(f(n))$ follows by the composition identity.

Unfortunately, the answer to the Question \ref{klaus} is negative. In the next section we prove a general result which holds for polynomials with coefficients in an arbitrary field $K$. Our result shows that one cannot prove the conjecture for cubic polynomials $f(x)$ by using the composition identity in Question \ref{klaus}. We also refer to \cite{dubmul}, where a certain composition identity was used to investigate multiplicative dependence of integer values of quadratic integer polynomials and \cite{dubdrun} for further results in this direction.

\section{Main Result}
The main result of this paper is the following theorem:
\begin{theorem}\label{quadr} Let $m \geq 2$ be an integer not divisible by the characteristic of the field $K$. Suppose that $f(x) \in K[x]$ is non constant and separable, and the polynomial $g(x)$, $\deg g \geq 2$, has a non-zero derivative. Then the equation
\[f(g(x)) = f(x)h^m(x)\]
\noindent holds if and only if:
\[I) \; \: \qquad f(x) = ax+b, \qquad a, b \in K, a \ne 0, \qquad g(x) = \left(x+\frac{b}{a}\right)h^m(x)-\frac{b}{a}\]
\noindent or
\[II) \qquad f(x) = ax^2 + bx +c, \qquad a, b, c \in K, \quad a \ne 0, \quad \qquad  m = 2, \quad \qquad\]
\noindent with
\[\qquad g(x) = \frac{1}{2a}\left(\pm T_n\left( \frac{2ax+b}{\sqrt{D}}\right)\sqrt{D}  -b\right), \qquad h(x) = \pm U_{n-1}\left(\frac{2ax+b}{\sqrt{D}}\right),\]

\noindent where $T_n(x)$, $U_n(x)$ are Chebyshev polynomials of the first and second kind, respectively, $D = b^2-4ac$ is the discriminant of $f(x)$.
\end{theorem}

We remark that the condition on the separability of $f(x)$ cannot be weakened in Theorem \ref{quadr} which can be seen by taking $f(x)=g(x) = x(x-1)^m$ in $\Q[x]$. The requirement that $g(x)$ has a non-zero derivative for fields $K$ of characteristic $p \ne 0$ also cannot be weakened. Indeed, consider the simple example given by $f(x)=x^d-1$, $g(x) = x^{p^l}$ in $\mathbb{F}_p[x]$. Also, if the characteristic $p$ divides the exponent $m \ne 0$ in the equation $f(g(x))=f(x)h^m(x)$, then one can write $h^m(x) = h_1^{m/p}(x^p) = h_2^{m/p}(x)$, where $h_2(x)$ is a polynomial with coefficients in $K$.

Recall that for the field $K$ of characteristic not equal to $2$, the \emph{Chebyshev polynomials} $T_n(x) \in K[x]$ of the \emph{first kind} are defined by the linear recurrence of order two:
\begin{equation}\label{defT}
T_0(x) = 1, \quad T_1(x) = x,  \qquad T_{n+2}(x) = 2xT_{n+1}(x)-T_n(x).
\end{equation}
In the similar way, the Chebyshev polynomials of the \emph{second kind} $U_n(x) \in K[x]$ are defined by the recurrence
\begin{equation}\label{defU}
U_0(x) = 1, \quad U_1(x) = 2x, \qquad U_{n+2}(x) = 2xU_{n+1}(x)-U_n(x).
\end{equation}

Polynomials $T_n(x)$ and $U_n(x)$ contain only even powers of $x$ for even $n$, odd powers of $x$ for odd $n$. Thus, the coefficients of $g(x)$ and $h(x)$ in Theorem \ref{quadr}, (II) lie in $K$ if $n$ is odd and in $K(\sqrt{D})$ if $n$ is even. Chebyshev polynomials have many other remarkable properties, see, for instance, \cite{rivl}. They play a key role in the theorems of Ritt for decompositions of polynomials \cite{schinz}. In addition, Chebyshev polynomials are related to permutation polynomials over finite fields called Dickson polynomials \cite{lidl}. In our proof, the following property of Chebyshev polynomials will be useful:

\begin{proposition}\label{pell1} Suppose that the characteristic of the field $ K$ is not equal to $2$. Then all solutions of the Pell equation
\[ P^2(x) - (x^2-1)Q^2(x) =1\]
in the ring $K[x]$ are given by \[ P(x) = \pm T_n(x), \qquad Q(x) =  \pm U_{n-1}(x),\]
where $T_n(x)$ and $U_n(x)$ are Chebyshev polynomials of the first and second kind, respectively.
\end{proposition}

The equation wich appears in Proposition \ref{pell1} is a special case of a general polynomial Pell equation $P(x)^2-D(x)Q^2(x) = 1$.  Solutions to general Pell equations in polynomials over complex number field $K=\C$ were investigated by Pastor \cite{past}. Dubickas and Steuding \cite{dubpell} gave an elementary algebraic proof for arbitrary field $K$. The proof of  Proposition \ref{pell1} can be found in \cite{dubpell}. Alternative proofs (in the case $K=\C$) are given in \cite{barb} and \cite{past}.

\section{Proof of Theorem \ref{quadr}}

\begin{proof}
Set $d = \deg f$. Let $a \in K$ and $b \in K$ be the leading coefficients of polynomials $f(x)$ and $g(x)$, respectively, $ab \ne 0$. Suppose that $L$ is the field extension of $K$ generated by the roots of the polynomials $f(x)$, $x^m-1$ and $x^m-b$. Then \begin{equation}\label{skaid}f(x) = a \prod_{\al \in V(f)} (x-\al).\end{equation} Here $V(f) \subset L$ denotes the set of the roots of the polynomial $f(x)$. The composition equation $f(g(x))=f(x)h^m(x)$ factors  in $L[x]$ into
\begin{equation}\label{sides2}
a\prod_{\al \in V(f)}(g(x) - \al) = a \prod_{\al \in V(f)}(x-\al) h^m(x),
\end{equation} and one can cancel $a$ on both sides. Observe that distinct factors $g(x)  - \al$ on the left hand side of \eqref{sides2} are relatively prime in $L[x]$ since their difference is a non-zero constant. We claim that at most one factor $g(x) - \al$ may be relatively prime with $f(x)$ if $m \geq 2$ and the characteristic of $K$ does not divide $m$. Indeed, suppose that $g(x) - \be$, $\be \in V(f)$, $\be \ne \al$ is another such factor. Then both $g(x) - \al$ and $g(x) - \be$ divide $h^m(x)$, so $g(x) - \al$ and $g(x) - \be$  must be the $m$-th powers of some polynomials $u(x)$ and $v(x)$ in $L[x]$ which divide $h(x)$, say, $g(x) - \al= u^m(x)$ and $g(x) -\be = v(x)^m$. (Note that $u(x)$ and $v(x)$ belong to  $L[x]$ since the field $L$ contains all roots of $f(x)$ and the $m$--th roots of the leading coefficient $b$ of the polynomial $g(x)$). Then $u(x)^m - v(x)^m = \be - \al$ is a non-zero constant polynomial. On the other hand,
\[
u^m(x) - v^m(x) = \prod_{j = 0}^{m-1}(u(x)-\zeta^j v(x)),
\] where $\zeta$ is a primitive $m$--th root of unity in $L$ and at least one of polynomials $u(x)-\zeta v(x)$ has degree greater than or equal to one which is impossible.

\noindent Now, suppose that $V(f) = \{\al_1, \al_2, \dots, \al_d\} .$ Let $V_j$ be the set containing all distinct common  roots of the polynomial $g(x) - \al_j$ and the polynomial $f(x)$,
\[
V_j := V(g(x)-\al_j) \cap V(f).
\] Then $g(x) -\al_j = f_j(x)u_j(x)$, where $u_j(x) \in L[x]$ and \[f_j(x) := \prod_{\al \in V_j} (x -\al).\] Note that $f_j(x)$ are all separable and and coprime in $L[x]$. Since $f(x)$ is also separable, the equation  \eqref{sides2} implies
\begin{equation}\label{fact2}
a\prod_{j=1}^d f_j(x) = f(x) \qquad \text{ and consequently, } \qquad \prod_{j=1}^d u_j(x) = h^m(x).
\end{equation}
The polynomials $u_j(x)$ are relatively prime, thus $u_j(x) = h_j^m(x)$, $j=1$,$\dots$, $d$, for some polynomials $h_j(x) \in L[x]$ whose product is equal to $h(x)$ in \eqref{fact2}.
Let $n_j := \deg f_j$,  for $j=1,  \dots, d$. Without loss of generality, assume that $n_1 \leq n_2 \leq \dots \leq n_d$. Then $n_1 \geq 0$. Observe that $n_2 \geq 1$ if $n_1 = 0$, since no two factors $g(x)-\al_j$ can be coprime with $f(x)$, as noted above. The first identity in \eqref{fact2} gives
\begin{equation}\label{laipsn} n_1 + n_2 + \dots + n_d = \deg f = d.
\end{equation} Since $g(x)= f_j(x)h_j(x)^m+\al_j$, one also has $\deg g \equiv n_j \pmod{m}$.  We now consider two cases for $\deg{g}$ modulo $m$.

\medskip
\noindent \emph{Case 1).} Assume that $\deg g \equiv 0 \pmod{m}$. Then $n_j \geq m$ for $j \geq 2$, hence
\begin{equation}\label{nel1}
d \geq m(d-1)
\end{equation} by \eqref{laipsn}. Since $m \geq 2$,  one has $d \geq 2d -2$ which is possible for $d = 1$ or $d = 2$ only. Suppose that $d=2$. Then one also has $m \leq 2$ by \eqref{nel1}.
\medskip

\noindent \emph{Case 2).} Assume that $\deg g \not\equiv 0 \pmod{m}$. Then $n_1 = \dots = n_d =1$ by \eqref{laipsn}. Let $\deg g = sm + 1$, where $s: = \deg h_j \geq 1$ for $1 \leq j \leq d$. Since $h_j^m(x) \mid g(x) - \al_j$, the polynomials $h_j^{m-1}(x)$ are (relatively prime) factors of the derivative $g'(x)$. By conditions of Theorem, $g'(x)$ is a non-zero polynomial, hence \[ms \geq \deg g' \geq \deg h_1^{m-1} + \dots + \deg h_d^{m-1} = d(m-1)s\] and, consequently,
\begin{equation}\label{nel2}
m \geq d(m-1).
\end{equation}
Then $d \leq m/(m-1) \leq 2$. Suppose $d=2$. Then, in addition, \eqref{nel2} gives $m \leq 2$.

\medskip
\noindent Thus it remains to consider the cases $d = 1$ and $d = 2$. In the first case, the polynomial $f(x)$ is linear, thus $f(x) = ax+b$ with $a, b \in K$, $a \ne 0$. The equation $f(g(x))=f(x)h^m(x)$ is equivalent to
\[ ag(x)+b = (ax+b)h^m(x),\]
so one simplification solves $g(x)$ and this completes the proof in the case $d=1$.
\noindent Suppose $d=2$. Then $f(x) = ax^2+bx+c$ with $a, b, c \in K$, $a \ne 0$. Let $D = b^2-4ac$, $D \ne 0$ since $f(x)$ is separable. One also has $m=2$ by the conditions of Theorem \ref{quadr} and the degree inequalities in the two cases above. Hence, it suffices to find the polynomials $g(x)$ and $h(x)$ in the equation $f(g(x))=f(x)h^2(x)$. Since the characteristic of the field $K$ is not equal to $2$ by the conditions of Theorem \ref{quadr}, the linear change of variables $x \to x(t)$ defined by \[x = \frac{t\sqrt{D} -b}{2a}\] transforms the polynomial $f(x)$ into \[f(x) = \frac{D}{4a}F(t),\]
where $F(t)=t^2-1$. Set \[G(t) :=\frac{1}{\sqrt{D}} \left(2ag\left(\frac{t\sqrt{D}-b}{2a}\right)+b\right), \qquad H(t) :=h\left(\frac{t\sqrt{D}-b}{2a}\right).\]
By straightforward substitution, one easily checks that the map $x \to x(t)$ transforms the composition equation $f(g(x))=f(x)h^2(x)$ into $D/4a F(G(t))=D/4aF(t)H^2(t)$. Canceling the factor $D/4a$ on both sides, one obtains
\[F(G(t))=F(t)H^2(t),\]
or, equivalently,
\[G^2(t) - (t^2-1)H^2(t) = 1.\]
By Proposition \ref{pell1} all the solutions to this equation are given by the formulas $G(t) = \pm T_n(t)$, $H(t) = \pm U_{n-1}(t)$, where $T_n(t)$ and $U_n(t)$ are Chebyshev polynomials of the first and second kind, respectively. Application of the inverse map $t\to t(x)$ now yields the result.

\end{proof}

\section{Rational and integer examples}

Let $f(x) = ax^2 + bx+c$ be a quadratic polynomial with rational coefficients. For $n = 3$ in Theorem \ref{quadr}, one has $T_3(x) = 4x^3-3x$ and $U_2(x) = 4x^2-1$. Then $f(g(x))=f(x)h^{2}(x)$ holds by Theorem \ref{quadr} for
\begin{equation}\label{lyg}
\begin{aligned}
g(x) &= (16a^2x^3+24abx^2+(9b^2+12ac)x+8bc)/D,\\
h(x) &= (16a^2x^2+16abx+3b^2+4ac)/D.
\end{aligned}
\end{equation}
Extend the definition of $\la$ function to the whole set of rationals $\Q$ by the complete multiplicativity of $\lambda$. Then, using the method outlined in Section \ref{intr}, one can prove easily the following analogue of Theorem $2$ in \cite{lui} for the sign changes of $\la$ function at rational points $f(r)$, $r \in \Q$, namely: either $\la(f(r))$ is constant for all rational numbers $r$ greater than the largest real root of $g(x)-x$ or it changes sign infinitely many often.

The question of finding all solutions of the composition equation in integer polynomials $f(x)$, $g(x)$, and $h(x)$ is closely related to the solution of the polynomial Pell equations in $\Z[x]$, see \cite{mac}, \cite{nath}, \cite{wy}. This does not seem to be easy. The examples of such polynomials are $f(x) = x^2 \pm 1$, $f(x) = x^2 \pm 2$, $f(x) = x^2 \pm 4$. Respective polynomials $g(x)$ and $h(x)$ with integer coefficients can be found using \eqref{lyg}. See Table \ref{tb} bellow.
\begin{table}[h]
\caption{ Examples of polynomials $f(x), g(x), h(x) \in \Z[x]$ in Theorem \ref{quadr}. }\label{tb}
\centering
\begin{tabular}[c]{| c | c | c |}
\hline
 $f(x)$ & $g(x)$ & $h(x)$\\
\hline
$x^2+1$ & $4x^3+3x$ & $4x^2+1$\\
$x^2-1$ & $4x^3-3x$ & $4x^2-1$\\
$x^2+2$ & $2x^3+3x$ & $2x^2+1$\\
$x^2-2$ & $2x^3-3x$ & $2x^2-1$\\
$x^2+4$ & $x^3+3x$ & $x^2+1$\\
$x^2-4$ & $x^3-3x$ & $x^2-1$\\
\hline
\end{tabular}
\end{table}

\end{document}